\documentclass[12pt]{amsart}
\usepackage{amsmath, amsthm, amscd, amsfonts}

\newtheorem{theorem}{Theorem}[section]
\newtheorem{lemma}[theorem]{Lemma}
\newtheorem{proposition}[theorem]{Proposition}
\newtheorem{corollary}[theorem]{Corollary}
\theoremstyle{definition}
\newtheorem{definition}[theorem]{Definition}
\newtheorem{example}[theorem]{Example}

\theoremstyle{remark}

\numberwithin{equation}{section}
\begin{document}

\title{Generalized Derivations on Modules}

\author[Gh. Abbaspour, M. S. Moslehian, A. Niknam]{Gh. Abbaspour, M. S. Moslehian and A. Niknam}

\address{Gholamreza Abbaspour Tabadkan: Department of Mathematics, Ferdowsi University, P.O. Box 1159, Mashhad 91775, Iran}
\email{tabadkan@math.um.ac.ir}

\address{Mohammad Sal Moslehian: Department of Mathematics, Ferdowsi University, P.O. Box 1159, Mashhad 91775, Iran}
\email{moslehian@ferdowsi.um.ac.ir}

\address{Assadollah Niknam: Department of Mathematics, Ferdowsi University, P.O. Box 1159, Mashhad 91775, Iran}
\email{niknam@math.um.ac.ir}



\begin{abstract}
Let $A$ be a Banach algebra and $M$ be a Banach right $A$-module.
A linear map $\delta : M\to M$ is called a generalized derivation
if there exists a derivation $d : A \to A$ such that
$$\delta(xa)=\delta(x)a + x d(a) \quad (a \in A ,x \in M).$$
In this paper, we associate a triangular Banach algebra ${\mathcal
T}$ to Banach $A$-module $M$ and investigate the relation between
generalized derivations on $M$ and derivations on ${\mathcal T}$.
In particular, we prove that the so-called generalized first
cohomology group of $M$ is isomorphic to the first cohomology
group of ${\mathcal T}$.
\end{abstract}
\maketitle

\section{Introduction.}

Recently, a number of analysts \cite{A-M-N, B-V, M-M} have studied
various generalized notions of derivations in the context of
Banach algebras. There are some applications in the other fields
of study \cite{H-L-S}. Such maps have been extensively studied in
pure algebra; cf. \cite{A-R, BRE, HVA}.

Let, throughout the paper, $A$ denote a Banach algebra (not
necessarily unital) and let $M$ be a Banach right $A$-module.

A linear mapping $d : A \to A$ is called a derivation if $d(ab)=
d(a)b + ad(b)\quad(a, b\in A)$. If $a\in A$ and we define $d_{a}$
by $d_{a}(x)=ax-xa \quad (x\in A)$. Then $d_{a}$ is a derivation
and such derivation is called inner.

A linear mapping $\delta : M\to M$ is called a generalized
derivation if there exists a derivation $d : A \to A$ such that
$\delta(xa)=\delta(x)a + xd(a)\quad(x\in M, a\in A)$. For
convenience, we say that such a generalized derivation $\delta$
is a $d$-derivation. In general, the derivation $d : A\to A$ is
not unique and it may happen that $\delta$ (resp. $d$) is bounded
but $d$ (resp. $\delta$) is not bounded. For instance, assume that
the action of $A$ on $M$ is trivial, i.e $MA=\{0\}$. Then every
linear mapping $\delta : M\to M$ is a $d$-derivation for each
derivation $d$ on $A$.

Our notion is a generalization of both concepts of a generalized
derivation (cf. \cite{BRE, HVA}) and of a multiplier (cf.
\cite{DAL}) on an algebra (see also \cite{MOS2}). For seeing this,
regard the algebra as a module over itself. The authors in
\cite{A-M-N} investigated the generalized derivations on Hilbert
$C^*$-modules and showed that these maps may appear as the
infinitesimal generators of dynamical systems.

\begin{example} Let $M$ be a right Hilbert $C^{*}$-module over a $C^*$-algebra $A$ of compact operators acting on a
Hilbert space (see \cite{LAN} for more details on Hilbert
$C^*$-modules). By Theorem 4 of \cite{B-G}, $M$ has an
orthonormal basis so that each element $x$ of $M$ can be
expressed as
$x=\displaystyle{\sum_{\lambda}}v_{\lambda}<v_{\lambda},x>$. If
$d$ is a derivation on $A$, then the mapping $\delta : M\to M$
defined by $\delta(x)=\displaystyle{\sum
_{\lambda}}v_{\lambda}d(<v_{\lambda},x>)$ is a $d$-derivation
since
\begin{eqnarray*}
\delta(xa)&=&\delta\left(\displaystyle{\sum_{\lambda}}v_{\lambda}<v_{\lambda},xa>\right)\\
&=&\displaystyle{\sum_{\lambda}}v_{\lambda}d(<v_{\lambda},x>a)\\
&=&\displaystyle{\sum_{\lambda}}v_{\lambda}d(<v_{\lambda},x>)a+\displaystyle{\sum_{\lambda}}v_{\lambda}<v_{\lambda},x>d(a)\\
&=&\delta(x)a + xd(a).
\end{eqnarray*}
\end{example}

The set ${\mathcal B}(M)$ of all bounded module maps on $M$ is a
Banach algebra and $M$ is a Banach ${\mathcal B}(M)-A$-bimodule
equipped with $T.x=T(x)\quad(x\in M, T\in {\mathcal B}(M))$, since
we have $T.(xa)=T(xa)=T(x)a=(T.x)a$ and $\|T.xa\|
\leq\|T\|\;\|x\|\;\|a\|$, for all $a \in A, x\in M, T\in {\mathcal
B}(M)$.

We call $\delta : M\to M$ a generalized inner derivation if there
exist $a\in A$ and $T\in {\mathcal B}(M)$ such that
$\delta(x)=T.x - xa = T(x)- xa$. Mathieu in \cite{MAT} called a
map $\delta : A\to A$ a generalized inner derivation if
$\delta(x)=bx-xa$ for some $a, b\in A$. If we consider $A$ as a
right $A$-module in a natural way, and take $T(x)=bx$, then our
definition covers the notion of Mathieu.

In this paper we deal with the derivations on the triangular
Banach algebras of the form ${\mathcal T}=\left(\begin{array}{cc}
{\mathcal B}(M) & M\\ 0 & A \end{array}\right)$. Such algebras
were introduced by Forrest and Marcoux \cite{F-M1} that in turn
are motivated by work of Gilfeather and Smith in \cite{G-S}
(these algebras have been also investigated by Y. Zhang who
called them module extension Banach algebras \cite{ZHA}). Among
some facts on generalized derivations, we investigate the
relation between generalized derivations on $M$ and derivations
on ${\mathcal T}$. In particular, we show that the generalized
first cohomology group of $M$ is isomorphic to the first
cohomology group of ${\mathcal T}$.

\section{Main Results.}

If we consider $A$ as an $A$-module in a natural way then we have
the following lemma about generalized derivations on $A$.

\begin{lemma}
A linear mapping $\delta : A \to A$ is a generalized derivation if
and only if there exist a derivation $d : A\to A$ and a module map
$\varphi :  A\to A$ such that $\delta=d+\varphi$.\end{lemma}
\begin{proof} Suppose $\delta$ be a generalized derivation on $A$,
then there exists a derivation $d$ on $A$ such that $\delta$ is a
$d$-derivation. Put $\varphi=\delta-d$. Then for each $a, x\in A$
we have
$$\varphi(xa)=\delta(xa)-d(xa)=\delta(x)a+xd(a)-(d(x)a+xd(a))=(\delta(x)-d(x))a=\varphi(x)a$$
Thus $\varphi$ is a module map and $\delta=d+\varphi$.

Conversely, let $d$ be a derivation on $A$, $\varphi$ be a module
map on $A$ and put $\delta=d+\varphi$. Then clearly $\delta$ is a
linear map and
$$\delta(xa)=d(xa)+\varphi(xa)=d(x)a+xd(a)+\varphi(x)a=(d(x)+\varphi(x))a+xd(a)=\delta(x)a+xd(a)$$
for all $a, x\in A$. Therefore $\delta$ is a
$d$-derivation.\end{proof}

The next two results concern the boundedness of a generalized derivation.

\begin{theorem} Let $A$ have a bounded left approximate identity $\{e_{\alpha}\}_{\alpha\in I}$
and let $\delta$ be a $d$-derivation on $A$. Then $\delta$ is
bounded if and only if $d$ is bounded. \end{theorem}
\begin{proof} First we show that every module map on $A$ is
bounded. Suppose that $\varphi$ is a module map on $A$ and let
$\{a_{n}\}$ is a sequence in $A$ converging to zero in the norm
topology. By a consequence of Cohen Factorization Theorem (see
Corollary 11.12 of \cite{B-D}) there exist a sequence $\{b_{n}\}$
and an element $c$ in $A$ such that $b_{n}\to 0$ and
$a_{n}=cb_{n},\quad (n\in {\mathbb N})$. Then
$\varphi(a_{n})=\varphi(cb_{n})=\varphi(c)b_{n}\to 0$. Thus by the
closed graph theorem, $\varphi$ is bounded.

Now let $\delta$ be a $d$-derivation. By Lemma 2.1,
$\delta=d+\varphi$ for some module map $\varphi$ on $A$.
Therefore $\delta$ is bounded if and only if $d$ is
bounded.\end{proof}

\begin{corollary} Every generalized derivation on a $C^{*}$-algebra is bounded.\end{corollary}
\begin{proof} Every derivation on a $C^*$-algebra is automatically continuous; cf. \cite{J-S}.\end{proof}

Let $\varphi : A\to A$ be a homomorphism (algebra morphism). A
linear mapping $T : M\to M$ is called a $\varphi$-morphism if
$T(xa)=T(x)\varphi(a)\quad(a\in A, x\in M)$. If $\varphi$ is a
isomorphism and T is a bijective mapping then we say T to be a
$\varphi$-isomorphism. An $id_{A}$-morphism is a module map (module
morphism). Here $id_{A}$ denotes the identity operator on $A$.

\begin{proposition} Suppose $\delta$ is a bounded $d$-derivation on $M$ and $d$ is bounded. Then $T=\exp(\delta)$ is
a bi-continuous $\exp(d)$-isomorphism.\end{proposition}
\begin{proof} Using induction one can easily show that
$$\delta^{(n)}(xa)=\sum_{r=0}^{n}(_{r}^{n})\delta^{(n-r)}(x)d^{(r)}(a).$$

For each $a\in A, x\in M$ we have
\begin{eqnarray*}
T(xa)&=&\exp(\delta)(xa)\\
&=&\sum_{n=0}^{\infty}\frac{1}{n!}\delta^{(n)}(xa)\\
&=&\sum_{n=0}^{\infty}\frac{1}{n!}\sum_{r=0}^{n}(_{r}^{n})\delta^{(n-r)}(x)d^{(r)}(a)\\
&=&\sum_{n=0}^{\infty}\sum_{r=0}^{n}(\frac{1}{(n-r)!} \delta^{(n-r)}(x)(\frac{1}{r!}d^{(r)}(a))\\
&=&(\sum_{n=0}^{\infty}\frac{1}{n!}\delta^{(n)}(x))(\sum_{n=0}^{\infty}\frac{1}{n!}d^{(n)}(a))\\
&=&\exp(\delta)(x)\exp(d)(a)
\end{eqnarray*}
The operators $\exp(\delta), \exp(d)$ are invertible in the
Banach algebras of bounded operators on $M$ and $A$,
respectively. Hence $T$ is an $\exp(d)$-isomorphism.\end{proof}

\begin{proposition}
Let $\delta$ be a bounded generalized derivation on $M$. Then
$\delta$ is a generalized inner derivation if and only if there
exists an inner derivation $d_{a}$ on $A$ such that $\delta$ is
$d_{a}$-derivation. \end{proposition}

\begin{proof} Let $\delta$ be a generalized inner derivation. Then there exist $a\in A$ and $T\in
{\mathcal B}(M)$ such that $\delta(xa)=T(x)-xa\quad (x\in M)$. We
have $\delta(x)b+xd_{a}(b)=(T(x)-xa)b+xab-xba
=T(x)b-xba=T(xb)-(xb)a=\delta(xb) \quad (b \in A, x \in M)$.
Hence $\delta$ is a $d_{a}$-derivation.

Conversely, suppose $\delta$ is a $d_{a}$-derivation for some
$a\in A$. Define $T : M\to M$ by $T(x)=\delta(x)+xa$. Then T is
linear, bounded and
$T(xb)=\delta(xb)+(xb)a=(\delta(x)b+xd_{a}(b))+xba=\delta(x)b+xab-xba+xba=(\delta(x)+xa)b=T(x)b$.
It follows that $T\in {\mathcal B}(M)$ and
$\delta(x)=(\delta(x)+xa)-xa=T(x)-xa$. Therefore $\delta$ is a
generalized inner derivation.\end{proof}

The linear spaces of all bounded generalized derivations and
generalized inner derivations on $M$ are denoted by $GZ^{1}(M,M)$
and $GN^{1}(M,M)$, respectively. We call the quotient space
$GH^{1}(M,M)=GZ^{1}(M,M)/GN^{1}(M,M)$ the generalized first
cohomology group of $M$.

\begin{corollary}
$GH^{1}(M,M)=0$ whenever $H^{1}(A,A)=0$
\end{corollary}
\begin{proof} Let $\delta : M\to M$ be a generalized derivation.
Then there exists a derivation $d : A \to A$ such that $\delta$ is
a $d$-derivation. Due to $H^{1}(A,A)=0$, we deduce that $d$ is
inner and, by Proposition 2.5, so is $\delta$. Hence
$GH^{1}(M,M)=0$.
\end{proof}

Using some ideas of \cite{F-M1, MOS1}, we give the following
notion:

\begin{definition}{\rm ${\mathcal T}=\{\left(\begin{array}{cc}T & x \\0 & a \end{array}\right);
T\in {\mathcal B}(M), x\in M, a\in A \}$ equipped with the usual
$2\times 2$ matrix addition and formal multiplication and  with
the norm $\|\left(\begin{array}{cc} T & x \\ 0 & a
\\\end{array}\right)\|=\|T\|+\|x\|+\|a\|$ is a Banach algebra. We
call this algebra the triangular Banach algebra associated to M.}
\end{definition}

The following two theorems give some interesting relations
between generalized derivations on $M$ and derivations on
${\mathcal T}$.

Let $\delta$ be a bounded $d$-derivation on $M$. We define
$\Delta_{\delta} : {\mathcal B}(M)\to {\mathcal B}(M)$ by
$\Delta_{\delta}(T)=\delta T-T\delta$. Then $\Delta_{\delta}$ is
clearly a derivation on ${\mathcal B}(M)$.

\begin{theorem}
Let $\delta$ be a bounded $d$-derivation on $M$ and $d$ be
bounded. Then the map $D^{\delta} : {\mathcal T} \to {\mathcal T}$
defined by $D^{\delta}\left(\begin{array}{cc} T & x \\ 0 & a
\end{array}\right)=\left(
\begin{array}{cc} \Delta_{\delta}(T) & \delta(x) \\  0 & d(a)\end{array}\right)$ is a bounded derivation on ${\mathcal T}$.
Also $\delta$ is a generalized inner derivation if and only if $D^{\delta}$ is an inner derivation. \end{theorem}

\begin{proof} It is clear that $D^{\delta}$ is linear. For any $T_{1},T_{2}\in {\mathcal B} (M), x_{1},x_{2}\in M, a_{1},a_{2}\in A$ we have
\begin{eqnarray*}
&&D^{\delta}(\left(\begin{array}{cc}T_{1} & x_{1} \\ 0 &
a_{1}\end{array}\right)\left(\begin{array}{cc} T_{2} & x_{2} \\ 0 & a_{2} \end{array}\right))=
D^{\delta}\left(\begin{array}{cc} T_{1}T_{2} & T_{1}.x_{2}+x_{1}a_{2} \\
0 & a_{1}a_{2} \end{array}\right)\\
&=&\left(\begin{array}{cc} \Delta_{\delta}(T_{1}T_{2}) &
\delta(T_{1}.x_{2}+x_{1}a_{2}) \\  0 & d(a_{1}a_{2})\end{array}
\right )\\
&=&\left(\begin{array}{cc} \Delta_{\delta}(T_{1}T_{2}) &
\delta(T_{1}(x_{2}))+ \delta(x_{1})a_{2}+x_{1}d(a_{2})\\ 0 &
a_{1}d(a_{2})+d(a_{1})a_{2}\end{array}\right )\\
&=&\left(\begin{array}{cc}
T_{1}\Delta_{\delta}(T_{2})+\Delta_{\delta}(T_{1})T_{2} &
T_{1}.\delta(x_{2})+x_{1}d(a_{2})+(\delta
T_{1}-T_{1}\delta)(x_{2})+\delta(x_{1})a_{2} \\ 0 & a_{1}d(a_{2})+d(a_{1})a_{2}\end{array}\right )\\
&=&\left(\begin{array}{cc} T_{1} & x_{1} \\  0 & a_{1}
\end{array}\right)\left(\begin{array}{cc} \Delta_{\delta}(T_{2}) & \delta(x_{2})
\\ 0 & d(a_{2})\end{array}\right) + \left(\begin{array}{cc} \Delta_{\delta}
(T_{1}) & \delta(x_{1}) \\ 0 & d(a_{1}) \end{array}\right)\left(\begin{array}{cc} T_{2} & x_{2} \\ 0 & a_{2}\end{array}\right)\\
&=&\left(\begin{array}{cc} T_{1} & x_{1} \\ 0 & a_{1}
\end{array}\right)D^{\delta}\left(\begin{array}{cc} T_{2} & x_{2}
\\ 0 & a_{2}\end{array}\right)+D^{\delta}(\left(\begin{array}{cc}  T_{1} & x_{1} \\ 0 & a_{1}\end{array}\right))\left(
\begin{array}{cc} T_{2} & x_{2} \\ 0 & a_{2} \end{array}\right)
\end{eqnarray*}

Thus $D^{\delta}$ is a derivation on ${\mathcal T}$. Due to
$\|\left(\begin{array}{cc} \Delta_{\delta}(T) & \delta(x) \\
0 & d(a)
\end{array}\right)\|=\| \Delta_{\delta}(T)\|+\|\delta(x)\|+\|d(a)\|\leq
\max \{\|\Delta_{\delta}\|,\|\delta\|,\|d\|\}\|\left(\begin{array}{cc} T & x \\
0 & a \end{array}\right)\|$, we infer that $D^{\delta}$ is
bounded.

Now suppose that $\delta$ is a generalized inner derivation. Then
there exist $a\in A$ and $T\in {\mathcal B}(M)$ such that
$\delta(x)=T(x)-xa  \quad(x\in M)$. For all $S\in {\mathcal
B}(M), b\in A$ and $y\in M $ we have
\begin{eqnarray*}
D_{\left(\begin{array}{cc} T & 0
\\ 0 & a \end{array}\right)}\left(\begin{array}{cc}  S & y \\  0
& b \end{array}\right) &:=& \left(\begin{array}{cc} T & 0 \\  0 &
a \end{array} \right)\left(\begin{array}{cc} S & y \\ 0 & b
\end{array}\right)-\left(\begin{array}{cc} S & y \\ 0 & b\end{array}\right)\left(\begin{array}{cc} T & 0 \\ 0 & a
\end{array} \right)\\
&=&\left(\begin{array}{cc} TS-ST & T.y-ya \\ 0 & ab-ba
\end{array}\right)\\
&=& \left(\begin{array}{cc} \Delta_{\delta}(S) & \delta(y)
\\0 & d_{a}(b)\end{array}\right)\\
&=&D^{\delta}\left(\begin{array}{cc} S & y \\ 0 & b
\end{array}\right). \end{eqnarray*}

Hence $D^{\delta}=D_{\left(\begin{array}{cc} T & 0 \\ 0 & a
\end{array} \right)}$ and so $D^{\delta}$ is an inner derivation.

Conversely, let $\delta$ be a bounded $d$-derivation such that
the associated derivation $D^{\delta}$ be an inner derivation,
say $D^{\delta} = D_{\left(\begin{array}{cc} T_0 & x_0
\\ 0 & a_0 \end{array}\right)}$. Then for each $T \in {\mathcal B}(M), x \in
M,a \in A$ we have
\begin{eqnarray}\label{inner}
\left(\begin{array}{cc} \Delta_\delta(T) & \delta(x) \\ 0 & d(a)
\end{array} \right)
&=& D^{\delta}(\left(\begin{array}{cc} T & x \\ 0 & a \end{array} \right))\nonumber \\
&=& D_{\left(\begin{array}{cc} T_0 & x_0
\\ 0 & a_0 \end{array}\right)}\left(\begin{array}{cc}  T & x \\  0
& a \end{array}\right)\nonumber \\
&=& \left(\begin{array}{cc} T_0T-TT_0 & T_0(x)+x_0a-T(x_0)-xa_0 \\
0 & a_0a-aa_0 \end{array}
\right)\nonumber \\
&=&\left(\begin{array}{cc} T_0T-TT_0 & T_0(x)+x_0a-T(x_0)-xa_0
\\0 & d_{a_0}(a)
\end{array}\right)
\end{eqnarray}
Hence $d=d_{a_0}$ is inner. Putting $a=0$ and $T=0$ in
(\ref{inner}) we conclude that $\delta(x) = T_0(x) - xa_0 \quad
(x \in M)$. Hence $\delta$ is a generalized inner derivation.
\end{proof}

The converse of the above theorem is true in the unital case.

\begin{theorem} Let $A$ be unital and ${\mathcal T}$ be the triangular Banach algebra associated to a unital Banach
right $A$-module $M$. Assume that $D : {\mathcal T} \to {\mathcal
T}$ is a bounded derivation. Then there exist $m_{0}\in M$, a
bounded derivation $d : A\to A$ and a bounded $d$-derivation
$\delta : M\to M$ such that
$$D\left(\begin{array}{cc} T & x \\  0 & a \end{array}\right)=\left(\begin{array}{cc} \Delta_{\delta}(T)&
\delta(x)+m_{0}a-T.m_{0} \\
0 & d(a) \end{array}\right)$$ Moreover, $D$ is inner if and only
if $\delta$ is a generalized inner derivation.\end{theorem}

\begin{proof} We use some ideas of Proposition 2.1 of
\cite{F-M1}. By simple computation one can verify that

(i) $D\left(\begin{array}{cc} 0 & 0 \\ 0 & 1_{A}
\end{array}\right)=\left(\begin{array}{cc} 0 & m_{0} \\ 0 & 0
\end{array}\right)$ for some $m_{0}\in M$;

(ii) $D\left(\begin{array}{cc} 0 & 0 \\ 0 & a \end{array}\right)=\left(\begin{array}{cc} 0 & m_{0}a \\
0 & d(a) \end{array}\right)$ for some bounded derivation $d$ on
$A$;

(iii) $D\left(\begin{array}{cc} 0 & x \\ 0 & 0 \end{array}\right)=\left(\begin{array}{cc} 0 & \delta(x) \\
0 & 0 \end{array}\right)$ for some linear mapping $\delta$ on $M$;

(iv) $D\left(\begin{array}{cc} T & 0 \\ 0 & 0
\end{array}\right)=\left(\begin{array}{cc}
 \Delta_{\delta}(T) & -T.m_{0} \\ 0 & 0 \end{array}\right)$;

and finally $D\left(\begin{array}{cc} T & x \\  0 & a
\end{array}\right)=\left(\begin{array}{cc} \Delta_{\delta}(T) &
\delta(x)+m_{0}a-T.m_{0} \\ 0 & d(a)\end {array}\right)$.

We have
\begin{eqnarray*}
\left(\begin{array}{cc} 0 & \delta(xa) \\ 0 & 0
\end{array}\right)&=&D(\left(
\begin{array}{cc} 0 & xa \\ 0 & 0 \end{array}\right))=D(\left(\begin{array}{cc}
 0 & x \\ 0 & 0 \end{array}\right)\left(\begin{array}{cc} 0 & 0 \\ 0 & a \end{array}\right))\\
&=&\left( \begin{array}{cc} 0 & x \\ 0 & 0 \end{array}\right)D(\left(\begin{array}{cc} 0 & 0 \\
 0 & a \end{array}\right))+D(\left(\begin{array}{cc} 0 & x \\ 0 & 0 \end{array}\right))\left(
\begin{array}{cc} 0 & 0 \\ 0 & a \end{array}\right)\\
&=&\left(\begin{array}{cc} 0 & x \\ 0 & 0\end{array}\right)\left(\begin{array}{cc} 0& m_0a \\
0 & d(a) \end{array}\right)+\left(\begin{array}{cc} 0 & \delta(x) \\
0 & 0 \end{array}\right)\left( \begin{array}{cc}
0 & 0 \\ 0 & a \end{array}\right)\\
&=&\left(\begin{array}{cc} 0 & \delta(x)a+xd(a) \\ 0 & 0
\end{array}\right)
\end{eqnarray*}

Thus $\delta(xa)=\delta(x)a+xd(a)$ and so $\delta$ is a
$d$-derivation.

It is clear that $D$ is inner if and only if $d$ is inner and,
using Proposition 2.5, the latter holds if and only if $\delta$ is
a generalized inner derivation.\end{proof}

\begin{theorem} Let $A$ be a unital Banach algebra, $M$ be a unital Banach
right $A$-module and ${\mathcal T}=\left(\begin{array}{cc}
{\mathcal B}(M) & M\\ 0 & A \end{array}\right)$. Then
$H^{1}({\mathcal T},{\mathcal T})\cong GH^{1}(M,M)$\end{theorem}

\begin{proof} Let $\Psi : GZ^{1}(M,M)\to H^{1}({\mathcal T},{\mathcal T})$ be defined by
$$\Psi(\delta)=[D^{\delta}]$$ where $[D^{\delta}]$ represents the equivalence class of
$D^{\delta}$ in $H^{1}({\mathcal T},{\mathcal T})$. Clearly
$\Psi$ is linear. We shall show that $\Psi$ is surjective. To end
this, assume that $D$ is a bounded derivation on ${\mathcal T}$.
Let $\delta$, $d$, $\Delta_{\delta}$ and $m_{0}\in M$ be as in
the Theorem 2.9. Then
\begin{eqnarray*}
(D-D^{\delta})\left(\begin{array}{cc}
T & x \\ 0 & a \end{array}\right)&=&\left(\begin{array}{cc} \Delta_{\delta}(T) & \delta(x)+m_{0}a-T.m_{0}\\
0 & d(a) \end{array}\right) -\left(\begin{array}{cc}
\Delta_{\delta}(T) & \delta(x)\\ 0 & d(a)
\end{array}\right)\\
&=&\left(\begin{array}{cc} 0 & m_{0}a-T.m_{0}\\ 0 & 0
\end{array}\right)\\
&=&D_{\left(\begin{array}{cc} 0 & -m_{0} \\ 0 & 0
\end{array}\right)}\left(
\begin{array}{cc} T & x \\ 0 & a \end{array}\right).
\end{eqnarray*}
So $[D]=[D^{\delta}]=\Psi(\delta)$ and thus $\Psi$ is
surjective. Therefore $H^{1}({\mathcal T},{\mathcal T})\cong
GZ^{1}(M,M)/Ker(\Psi)$.

Note that  $\delta\in Ker(\Psi)$ if and only if $D^{\delta}$ is
inner derivation on ${\mathcal T}$. Hence
$Ker(\Psi)=GN^{1}(M,M)$, by Theorem 2.8. Thus $H^{1}({\mathcal
T},{\mathcal T})\cong GH^{1}(M,M)$.
\end{proof}

\begin{example} Suppose that $A$ is unital and $M=A$. Then ${\mathcal B}(A)=A$ and
so $GH^{1}(A,A)\cong H^{1}(\left(\begin{array}{cc}A&A\\0&A
\end{array}\right ), \left(\begin{array}{cc}A&A\\0&A
\end{array}\right ))=H^1(A,A)$, by Proposition 4.4 of
\cite{F-M2}. In particular, every generalized derivation on a
unital commutative semisimple Banach algebra \cite{S-W}, a unital
simple $C^*$-algebra \cite{SAK}, or a von Neumann algebra
\cite{KAD} is generalized inner.\end{example}

We have investigated the interrelation between generalized
derivations on a Banach algebra and its ordinary derivations. We
also studied generalized derivations on a Banach module in virtue
of derivations on its associated triangular Banach algebra. Thus,
we established a link between two interesting research areas:
Banach algebras and triangular algebras.

\textbf{Acknowledgment.} The authors sincerely thank the referee for
valuable suggestions and comments.

\end{document}